\documentclass[12pt]{amsart}
\usepackage{amssymb}

\newtheorem{theorem}{Theorem}[section]

\newtheorem{lemma}{Lemma}[section]
\newtheorem{cor}{Corollary}[section]

\title[Bounds for discrete tomography solutions]
{Bounds for \\ discrete tomography solutions}

\author[Birgit Van Dalen, Lajos Hajdu, Rob Tijdeman]{Birgit Van Dalen, Lajos Hajdu, Rob Tijdeman}

\thanks{Research supported in part by the Hungarian Academy of Sciences, by the OTKA grants K75566, K100339, and NK101680, and by the T\'AMOP 4.2.1/B-09/1/KONV-2010-0007 project. This project is implemented through the New Hungary Development Plan, co-financed by the European Social Fund and the European Regional Development Fund.}



\address{B.E. Van Dalen, R. Tijdeman\hfill\break
\indent Mathematisch Instituut, Leiden University\hfill\break
\indent Postbus 9512\hfill\break 
\indent 2300 RA Leiden\hfill\break 
\indent The Netherlands}
\email{bevandalen@gmail.com}
\email{tijdeman@math.leidenuniv.nl}

\address{L. Hajdu\hfill\break
\indent University of Debrecen, Institute of Mathematics\hfill\break
\indent and the Number Theory Research Group\hfill\break
\indent of the Hungarian Academy of Sciences\hfill\break
\indent Debrecen\hfill\break \indent P.O. 12. \hfill\break
\indent H-4010\hfill\break \indent Hungary}
\email{hajdul@math.unideb.hu}

\begin{document}

\begin{abstract} We consider the reconstruction of a function on a finite subset of $\mathbb{Z}^2$ if the line sums in certain directions are prescribed. The real solutions form a linear manifold, its integer solutions a grid. First we provide an explicit expression for the projection vector from the origin onto the linear solution manifold in the case of only row and column sums of a finite subset of $\mathbf{Z}^2$. Next we present a method to estimate the maximal distance between two binary solutions. Subsequently we deduce an upper bound for the distance from any given real solution to the nearest integer solution. This enables us to estimate the stability of solutions. Finally we generalize the first mentioned result to the torus case and to the continuous case.
\end{abstract}

\maketitle
 
\section{Introduction}

The basic problem of discrete tomography is to reconstruct a function $f: A \to B$ where $A$ is a finite subset of $\mathbb{Z}^l$ and $B$ a finite subset of $\mathbb{R}$, when the sums of the function values along the lines in a finite number of directions are given. In this paper we consider $l=2$. There is a vast literature on the special case $A= \{ (i,j) \in \mathbb{Z}^2 | 0 \leq i <m, 0 \leq j < n \}$, $B = \{0, 1\}$ where the problem is to find the function values from the given row and column sums. In 1957 Ryser \cite{ry} and Gale \cite{ga} independently derived necessary and sufficient conditions for the existence of a solution in this special case. Ryser also provided a polynomial time algorithm for finding such a solution. However, the problem is usually highly underdetermined and a large number of solutions may exist\ \cite{wz}. Therefore the quest is often to find a solution of a certain type. For some classes of highly-structured images, such as hv-convex polynominoes (the 1's in each row and column are contiguous) polynomial time reconstruction algorithms have been developed (see e.g.\ \cite{bdnp}, \cite{bddkn}, \cite{df}).  Woeginger \cite{wo} presented an overview of classes of polyminoes for which it is proved that they can be reconstructed in polynomial time or it is proved that reconstruction is NP-hard. Batenburg \cite{ba1} developed an evolutionary algorithm for finding the reconstruction which maximises an evolution function and showed that the algorithm can be successfully applied to a wide range of evolution functions. 

We consider solutions as vectors with the values of $f$ as entries. 
Hajdu and Tijdeman \cite{ht1} observed that the set of binary solutions is precisely the set of shortest vector solutions in the set of functions $f: A \to \mathbb{Z}$ with the given line sums, provided that such solutions exist. They also showed that the solutions $f: A \to \mathbb{Z}$ with the given line sums form a multidimensional grid on the linear manifold which consists of all the solutions $f: A \to \mathbb{R}$ with the given line sums. Moreover, they determined the dimension of this manifold and indicated how to find a set of generators of the grid. Later they used their analysis to develop an algorithm to actually construct solutions $f: A \to {\{0,1\}}$ in \cite{ht2}, whereafter Batenburg \cite{ba2}, \cite{ba3} constructed much faster algorithms.

For functions $f: \{0,1, \dots, m-1\} \times \{0,1, \dots, n-1\} \to \mathbb{R}$ and given row and column sums we deduce a new and explicit expression for the projection vector $\vec{f_0}$ from the origin onto the real linear solution manifold in Section 3. 
We do not know a similar expression for other sets of line sums, but show the result can be extended to the torus case ($A = (\mathbb{Z} / n\mathbb{Z})^2$)
and to the continuous case ($A = [0,m] \times [0,n] \subset \mathbb{R}^2$) in Sections 6 and 7, respectively. In these sections we answer questions posed by Joost Batenburg.

In many applications it  suffices to find a solution or almost-solution which is guaranteed to be similar to the original. If all the solutions are similar, then they will also be similar to the original and we say that the solution set is stable. Alpers, Gritzmann and Thorens \cite{agt} showed that already a small change in the data can lead to a dramatic change in the image. Research on uniqueness and stability of solutions  was carried out by Alpers et al. \cite{al}, \cite{ag}, \cite{ab} and Van Dalen \cite{dal0}, \cite{dal1}, \cite{dal2}, \cite{dal3} in case only row sums and column sums are given. Their estimates depend on only few parameters. More general situations were studied by Brunetti and Daurat \cite{bd} and by Gritzmann, Langfeld and Wiegelmann \cite{glw}. In Section 4 we use the distance estimates for the solutions to derive new stability results. They involve more parameters, but, at least in the given example, yield better results than the estimates obtained by Van Dalen. 

In Stolk \cite{st} a system of line sums is called compatible if and only if a real solution exists. (Then the projection vector $\vec{f_0}$ is an example of such a real solution.) Stolk showed that if the line sums are integers and they are compatible, then there exists an integer solution. In Section 5 we derive an upper bound for the Euclidean distance from $\vec{f_0}$ to the nearest integer solution.

The interest in discrete tomography arose from the study of atom positions in a crystal, but the developed theory has also applications in medical imaging and in nuclear science, see \cite{hk1, hk2}.
The results in the present paper are merely of theoretical interest, but applications of the method can be found in \cite{bfht} and they may help to estimate how many directions are needed to be sure that the solution is unique so that one is certain to have found the original configuration.

\section{Notation and general results}

We use the same settings as in \cite{ht1}. Let $a,b \in \mathbb{Z}$ with gcd$(a,b) = 1$ and $a \geq 0$. We call $(a,b)$ a direction. Put
 \[ f_{a,b}(x,y) = \left\{ \begin{array} {ll}
 x^ay^b-1, & {\rm if}~a>0, b>0, \\
 x^ay^{-b}, & {\rm if}~a>0,b<0, \\
 x-1,  & {\rm if} ~ a=1, b=0, \\
 y-1, & {\rm if}~ a=0, b=1.
\end{array} \right\} \]
By lines with direction $(a,b)$ we mean lines of the form $ay=bx+t~(t \in \mathbb{Z})$ in the $(x,y)$-plane.

Let $m$ and $n$ be positive integers and
$$A:=\{(i,j)\in{\mathbb Z}^2:0\leq i<m,0\leq j<n\}.$$ Let $R$ be an integral domain such that $R[x,y]$ is a unique factorization domain. 
(We apply the results from this section both with $R= \mathbb{R}$ and with $R= \mathbb{Z}$.)
If $g: A \to R$ is a function, then the line sum of $g$ along the line $ay=bx+t$ is defined as $\sum g(i,j)$ where the sum is extended over all $(i,j) \in A$ with $aj=bi+t$.

We often consider a function $v$ defined on $A$ as a vector with the $mn$ values of $v$ as entries. If we want to emphasize this, we write $\vec{v}$ instead of $v$.
We always assume that the coordinates of these vectors are arranged according to the elements of $A$ in lexicographical order. 
If $R \subseteq \mathbb{R}$, the length of $\vec{v}$ (or $v$) equals $|v| = | \vec{v}| = \sqrt {\sum_{(i,j) \in A} v(i,j)^2}$.

A set of directions $S=\{(a_d,b_d)\}_{d=1}^k$ is called valid for $A$ if $\sum_{d=1}^k a_d <m$ and $\sum_{d=1}^k |b_d| < n.$ 
Suppose $S = \{(a_d,b_d) \}_{d=1}^k$ is a valid set of directions for $A$. Write $M = \sum_{d=1}^k a_d, N = \sum_{d=1}^k |b_d|$.
Put $$F_S(x,y) = \prod_{d=1}^k f_{(a_d,b_d)}(x,y)~ {\rm and}~ F_{(u,v,S)}(x,y) = x^uy^vF_S(x,y)$$ for $0 \leq u < m-M, 0 \leq v < n-N$.
For these values of $u$ and $v$ define the functions $m_{(u,v,S)}: A \to R$ by 
$$ m_{(u,v,S)}(i,j) = {\rm coeff.~of~}(x^iy^j)~ {\rm in}~ F_{(u,v,S)}~~~{\rm for}~(i,j) \in A.$$
The $m_{(u,v,S)}$'s are called the switching elements corresponding to the direction set $S$.
By the bottom-left corner of the switching element $m_{(u,v,S)}$ we mean the lexicographically first element of $A$ for which the function value of $m_{(u,v,S)}$ is nonzero.
It is obvious that $F_{(u,v,S)} \in \mathbb{Z}[x,y]$.
As noticed on page 122 of \cite{ht1}, it follows from the above definitions that the function value of $F_{(u,v,S)}$ at its bottom-left corner is $\pm 1$.

We refer to the situation as described above as the general case, and to the case that $k=2$ and the directions are $(1,0),(0,1)$ as the simple case. 
In the latter case we only have row and column sums. 

The following result characterizes the structure of the solution set over the reals. 

\begin{lemma} [\cite{ht1}, Theorem 1] Let $S$ be a valid set of directions for $A$. 
Then any function $g: A \to R$ with zero line sums along the lines corresponding to $S$ can be uniquely written in the form
$$g = \sum_{u=0}^{m-1-M} \sum_{v=0}^{n-1-N} c_{uv}m_{(u,v,S)} $$
with $c_{uv} \in R$. On the other hand, every such function has zero line sums along the lines corresponding to $S$.
\end{lemma}

%\begin{proof} Theorem 1 of \cite{ht1}.
%\end{proof}

\noindent Let $A, S$ and $f: A \to \mathbb{R}$ be given. Then we denote by  $\vec{f_0}$ the vector which is orthogonal to the $(m-M)(n-N)$-dimensional linear subspace of the vectors $\vec{g}$ 
such that $g$ has zero line sums with respect to $S$ and for which $f_0$ has the same line sums as $f$ with respect to $S$.

%\begin{proof} The formula for $g$ in Theorem 1 of \cite{ht1} implies that the dimension is at most $(m-M)(n-N)$, the fact that the expression is unique that the dimension equals $(m-M)(n-N)$.
%\end{proof}

We shall also use the following consequence of Lemma 2.1.
\begin{lemma} [\cite{ht1}, Corollary 1] Let $K$ be the set of  the bottom-left corners of the switching elements of $A$ with respect to $S$. 
Then for any $f: A \to R$ and for any prescribed values for the elements of $K$, 
there exists a unique $g: A \to \mathbb{Z}$ having the prescribed integer values at the elements of $K$ and having the same line sums as $f$ along the lines corresponding to $S$.
\end{lemma}

%\begin{proof} Corollary 1 of \cite{ht1}.
%\end{proof}

\noindent It further follows from Lemma 2.1 that there are exactly $MN - \sum_{d=1}^k a_d|b_d|$ linearly independent homogeneous linear dependencies among the line sums.

Let $A, S$ and the line sums along the lines with respect to $S$ be given. Checking whether there is a real solution is simple. It suffices to solve the corresponding set of linear equations and to see whether it admits a solution. If the line sums are integers and the system is solvable, then obviously $\vec{f_0}$ is a rational solution. However, it is not obvious that there is an integer solution. The existence of such solutions is guaranteed by a result of Stolk (cf. p. 20 of \cite{st}). 

\begin{lemma} [Stolk, \cite{st}, Corollary 3.2.10]
\label{exist}
Let $A, S$ be given. Suppose all the line sums into the directions of $S$ are integers and there exists a real function $g: A \to \mathbb{R}$ satisfying the line sums. 
Then there exists a function $f: A \to \mathbb{Z}$ satisfying the line sums.
\end{lemma}

\noindent For the convenience of the reader we sketch the proof.
\vskip.1cm

\noindent {\it Proof sketch.}
Denote the lines in the directions of $S$ by $L$ and the number of directions of $S$ by $d$.
The proof proceeds by omitting elements from $A$ which are known to have integer values whence it follows that the remaining elements have integer sums along the lines of $L$. 

If a line of $L$ contains only one element, then this element has an integer value. Therefore we can omit it from $A$ and the remaining elements have integer sums along the remaining lines of $S$.
Denote the remaining elements by $A$ and the lines of $S$ containing at least one element of this $A$ by $L$. If the number of remaining directions is less than $d$, we are finished. 

Otherwise every line of $L$ contains at least two elements of $A$. Consider the convex hull of $A$.
This polygon contains at least $2d$ line segments, two into each direction of $S$.
By construction the convex hull of the switching element also contains two line segments into each direction of $S$.
Moreover, by the gcd-conditions it is the smallest nontrivial polygon with this property.
Hence the convex hull of $A$ contains a switching element.
Consider the bottom-left corners of the switching elements contained in the convex hull of $A$ and take the lexicographically smallest among them.
According to Lemma 2.2 we can choose any integer value for this element. We omit this element from $A$. The remaining elements have integer line sums along $S$.

From now on we proceed by iteration until $A$ becomes the empty set. If there is a line with only one element of $A$ we omit this element. Otherwise we select the first bottom-left corner of a switching element contained in $A$, choose an integer value for it and omit it.

\section{Explicit expression for the shortest real solution \\ in the simple case}

In case of only row and column sums we give an explicit form of $f_0$. For this we simplify our notation. Let $c_i$ $(i=0,\dots,m-1)$ and $r_j$ $(j=0,\dots,n-1)$ denote the column sums and row sums, respectively. Further, write $D=\sum\limits_{j=0}^{n-1} r_j$. Note that we have $D=\sum\limits_{i=0}^{m-1} c_i$.

\begin{theorem}
\label{thm2}
For any $(i,j)\in A$ we have
$$
f_0(i,j)={\frac{c_i}{n}}+{\frac{r_j}{m}}-{\frac{D}{mn}}.
$$
\end{theorem}

\begin{proof} To prove the statement, we need to check two properties: \\ $f_0$ is a solution, and $\vec{f_0}$ is orthogonal to $H$.

We start with the first property. Obviously, for any $i=0,\dots,m-1$, we have
$$
\sum\limits_{j=0}^{n-1} f_0(i,j)=\sum\limits_{j=0}^{n-1}
\left({\frac{c_i}{n}}+{\frac{r_j}{m}}-{\frac{D}{mn}}\right)=
c_i-{\frac{D}{m}}+{\frac{1}{m}}\sum\limits_{j=0}^{n-1} r_j=c_i.
$$
Similarly, for any $j=0,\dots,n-1$,
$$
\sum\limits_{i=0}^{m-1} f_0(i,j)=\sum\limits_{i=0}^{m-1}
\left({\frac{c_i}{n}}+{\frac{r_j}{m}}-{\frac{D}{mn}}\right)=
r_j-{\frac{D}{n}}+{\frac{1}{n}}\sum\limits_{i=0}^{m-1} c_i=r_j,
$$
which confirms the first property.

To prove the second property we check orthogonality for arbitrary $\vec{h}\in H$. Since for any $\vec{h}\in H$ with the corresponding $h:\  A \to {\Bbb R}$ both
$$
\sum\limits_{j=0}^{n-1} h(i,j)=0\ \ \ (i=0,\dots,m-1)
$$
and
$$
\sum\limits_{i=0}^{m-1} h(i,j)=0 \ \ \ (j=0,\dots,n-1),
$$
we have
$$
\sum\limits_{i=0}^{m-1}\sum\limits_{j=0}^{n-1} f_0(i,j)h(i,j)=
\sum\limits_{i=0}^{m-1}\sum\limits_{j=0}^{n-1}
\left({\frac{c_i}{n}}+{\frac{r_j}{m}}-{\frac{D}{mn}}\right)
h(i,j)=
$$
$$
=\sum\limits_{i=0}^{m-1} \left({\frac{c_i}{n}}-{\frac{D}{mn}}\right)
\sum\limits_{j=0}^{n-1} h(i,j)+
\sum\limits_{j=0}^{n-1}{\frac{r_j}{m}}
\sum\limits_{i=0}^{m-1} h(i,j)=0,
$$
and the theorem follows.
\end{proof}

\noindent {\bf Remark.} For the proof of Theorem 3.1 it seems to be essential that every two lines into the same direction contain the same number of pixels, i.e. elements of $A$
%We do not know a similar result in the general case.
(cf. Sections 6 and 7).

\section{Binary solutions}

In this section we consider the general binary case $f: A \to \{0,1\}$. 
Let $D:= \sum_{(a,b) \in A} f(a,b)$.
Then $D$ equals the sum of all the line sums in an arbitrary direction $(a_d,b_d)$. 
We present a method to give a lower bound for the number of correct pixel values in an approximate solution and to give an upper bound for the number of entries where two solutions can be different.
Such results have been obtained in the simple case by Van Dalen \cite{dal2} by a completely different method. She calls a function $F_0: A \to \{0,1\}$ uniquely determined if there is no other function 
$A \to \{0,1\}$ having the same row and column sums. Let $\alpha(F)$ be half of the sum of the absolute differences between the row sums  of $F$ and the row sums of a uniquely determined function $F_0$ with the same column sums as $F$.
Then Van Dalen derived upper bounds for the number of places where two solutions $F_1$ and $F_2$ can differ in terms of $m, n, D=D(F_1)$ and $\alpha(F_1)$. We shall derive some other estimates and then compare these results with those in \cite{dal2}.

The following result, which is used in the proof of Theorem 2 of \cite{ht1}, shows that the Euclidean distance between $f_0$ and any binary solution is fixed.

\begin{theorem}
\label{thm3.1} 

For any solution $g: A \to \{0,1\}$ we have
$$
|\vec{g}-\vec{f_0}|=\sqrt{D-|\vec{f_0}|^2}.
$$
%In other words, the solutions $\vec{g}$ are on a hypersphere in the linear solution manifold, with center $\vec{f_0}$ and radius $\sqrt{D-|\vec{f_0}|^2}$.
\end{theorem}

\begin{proof}Observe that if $g$ is a binary solution then $|\vec{g}|=\sqrt{D}$. This means that such solutions are situated on a hypersphere with the origin as center, and of radius $\sqrt{D}$. 
According to Lemma 2.1 the solutions $\vec{g}$ are located on a linear manifold of dimension 
$(m-M)(n-N)$ orthogonal to $\vec{f_0}$. The intersection of this manifold and the hypersphere is a hypersphere (of the appropriate dimension) having $\vec{f_0}$ as center. By the theorem of Pythagoras we get that the radius of this hypersphere is $\sqrt{D-|\vec{f_0}|^2}$, and the theorem follows.
\end{proof}

Put $\langle x \rangle = \min (|x|,|1-x|)$ and $ E = \sum_{i=1}^m \sum_{j=1}^n  \langle f_0(i,j)  \rangle^2$.
Then the Euclidean distance between $ \vec{f_0} $ and the nearest integer vector with entries in $\{0,1\}$ is exactly $\sqrt{E}.$
Hence we have the following consequence of Theorem \ref{thm3.1}.

\begin{cor}If $ E + |\vec{f_0}|^2 > D$, then there is no binary solution. \\
 If $ E + |\vec{f_0}|^2 = D$, then the only solutions $ g: A \to \{ 0,1 \}$ are obtained by rounding each $f_0(i,j)$ to the integer $0$ or $1$ which is nearest to it.
\end{cor}

\noindent Note that the only entries where the rounding is not unique are those with value $1/2$. If such entries do not exist, the solution is unique.
\vskip.1cm

If $D- E - |f_0|^2 $ is positive, but not too large, we still may conclude that a certain fraction of the rounded values agrees with any solution $g:A \to \{0,1\}$. (In most cases we cannot tell which rounded values are correct and it may even be impossible to do so.) 
Suppose the rounded value $F(i,j) \in \{0,1\}$ of $f_0(i,j)$ is not the right value. If $x := f_0(i,j) \geq 1/2$, then, when replacing the value 1 by 0, in $E$ we have to replace $ \langle f_0(i,j) \rangle ^2= (1-x)^2$ by $x^2$.
%$$ x^2 = |x-1|^2 + (x^2 - |x-1|^2) =  \langle f_0(i,j) \rangle + (x - |x-1|)(x+|x-1|).$$
\noindent Hence the contribution increases by $2x-1 = 2 f_0(i,j)   -1$. Similarly, if $ f_0(i,j) \leq 1/2$, then the contribution to $E$ increases by $1-2x = 1-2 f_0(i,j) $.  

Order the values 
$ |2f_0(i,j) -1| $ in nondecreasing order, $b_1, b_2, \dots, b_{mn}$, say. According to Theorem \ref{thm3.1} $D - |f_0|^2$ equals $E$ plus the sum of the values $b_i$ which correspond to wrong values in $F$.
Let $s$ be the value with 
\begin {equation} \label{1}
 b_1 + \dots +b_s \leq D - E - |\vec{f_0}|^2 < b_1 + \dots + b_{s+1}.
 \end{equation}
Then at most $s$ pixels can have wrong values. Therefore at least $mn-s$ pixels have the right values. 
Similarly, let $t$ be the value with 
\begin{equation} \label{2}
b_1 + \dots +b_t \leq 2(D - E - |\vec{f_0}|^2) < b_1 + \dots + b_{t+1}.
\end{equation}
Then for any two solutions at most $t$ corresponding pairs of pixels can have different values. Therefore at least $mn-t$ such pixels have the same values. 
So we have derived the following result.

\begin{theorem}
\label{thm3.2}
{\rm (a)} Let $s$ be defined as above. For any solution $g: A \to \{0,1\}$ we have $g(i,j) = F(i,j)$ for at least $mn-s$ pairs $(i,j)$ $(i=1, \dots, m; j= 1, \dots, n).$ \\
{\rm (b)}  Let $t$ be defined as above. For any two solutions $g_1, g_2: A \to \{0,1\}$ we have $g_1(i,j) = g_2(i,j)$ for at least $mn-t$ pairs $(i,j)$ $(i=1, \dots, m; j= 1, \dots, n).$ 
\end{theorem}

This result can be slightly improved if there are line sums over the $F(i,j)$ which do not agree with the corresponding sum over the $f(i,j)$. Consider some direction with such line sums. 
We know that some values of $F$ have to be wrong and can therefore increase
the value of $E$ by securing that among $b_1,...,b_s$ and $b_1,...,b_t$ in (\ref{1})
and (\ref{2}), respectively, there are not too many representatives from some row.
This may yield smaller values of $s$ and $t$, hence better results.
\vskip.3cm

\noindent {\bf Example 4.1.} Let $m=6, n=5.$ Let the row sums be given by $5, 4, 3, 2, 1$ and the column sums by $4, 4, 3, 2, 1, 1$, respectively.
Then $D=15$ and, according to Theorem \ref{thm2}, $30 f_0$ is given by
\[ \begin{array}{rrrrrr}
34	&	34	&	28	&	22	&	16	&	16 \\
29     &	29	&	23	&	17	&	11	&	11 \\
24     &	24	&	18	&	12	&	6	&	6 \\
19	& 	19	&	13	&	7	&	1	&	1 \\
14	& 	14	&	8	&	2	&	-4	 &	-4 
\end{array} \]
\noindent Thus $|f_0|^2 = 166/15$. Further we get the following table $F$ by rounding the elements of $f_0$ to the nearest integer.
\[ \begin{array}{rrrrrr}
1	&	1	&	1	&	1	&	1	&	1	\\
1	&	1	&	1	&	1	&	0	&	0	\\
1	&	1	&	1	&	0	&	0	&	0	\\
1	&	1	&	0	&	0	&	0	&	0	\\
0	&	0	&	0	&	0	&	0	&	0	
\end{array} \] 
A calculation gives $E= 13/5$, hence $D - |\vec{f_0}|^2 -E = 4/3$. This yields $s=9, t=13$. Thus, by Theorem \ref{thm3.2}, every solution differs at most at 9 places from $F$ and any two solutions differ from each other at most at 13 places. In fact, the solutions below, due to Van Dalen \cite{dal2}, show that the actual numbers $s$ and $t$ can be as large as 8 and 10, respectively: 
\[ \begin{array}{rrrrrrrrrrrrrr}
1	&	1	&	1	&	1	&	0	&	1	& & & 1&1&1&1&1&0\\
1	&	1	&	1	&	0	&	1	&	0	& & & 1&1&1&1&0&0\\
1	&	1	&	0	&	1	&	0	&	0	& & & 1&1&1&0&0&0\\
1	&	0	&	1	&	0	&	0	&	0	& & & 1&1&0&0&0&0\\
0	&	1	&	0	&	0	&	0	&	0	& & & 0&0&0&0&0&1
\end{array} \] 

Theorem 2 of \cite{dal2} implies that two solutions can differ at most at 20 places. This estimate is worse than the 13 places obtained above. 
Van Dalen's estimate depends on few parameters, but holds only for the simple case.
Theorem 4.2 can be applied for any set of directions.

Since in the above example $F$ does not satisfy the required row sums, an improvement of the value of $E$ is possible. The column sums are correct, but the top row sum 6 is 1 too high and the bottom row sum 0 is 1 too low. 
For correct line sums only one of both 16's on the top row may be decreased and only one of both 14's on the bottom row may be increased.
It follows from the corresponding formula (\ref{1}) that every solution differs at most at 8 places from $F$.
The above examples show that 8 is the best possible.
The same reasoning applied to (\ref{2}) leads to an upper bound 11 for $t$.

The example can be generalized as follows. Let $n,q$ be positive integers
and set $m=(n+1)q$. Let the row sums be given by $c_j=n-1$ for $1 \leq j \leq q, c_{lq+j} = n-l$ for $1 \leq l \leq n, 1 \leq j \leq q$, and $c_j=1$
for $nq+1 \leq j \leq (n+1)q$. In \cite{dal2} Van Dalen showed that $t \leq 2q \sqrt{4nq(n+1)+1} -2q$ and that $t=2nq$ can be reached.
After an elaborate computation the argument in the present paper yields the estimates $s<2nq$ and $t<4nq-2q$. Therefore we have obtained an
improvement by roughly a factor $\sqrt{q}$ compared to the estimate for $t$ from \cite{dal2}.
Notice that the present bound for $t$ cannot be further improved by a factor 2.

\section{The function values are in $\mathbb{Z}$}
In this section we consider the case of line sums for functions \\ $f: A \to \mathbb{Z}$.
We derive  an upper bound for the solution nearest to some given function $h: A \to \mathbb{R}$ in the Euclidean sense.
The function $h$ can be considered as a prescribed model for the integer solution $f$ satisfying the line sums. 
By applying the result with $h=f_0$ we derive an upper bound for the shortest integer solution.

Let, as before,  $S = \{(a_d,b_d) \}_{d=1}^k $ be a set of directions. 
%Put $$M = \sum_{d=1}^k a_d, \\ N= \sum_{d=1}^k |b_d|, T = \sum_{d=1}^k a_d|b_d|.$$ We assume $m>M, n>N$.
Put again  $$F_S(x,y) = \prod_{d=1}^k f_{(a_d,b_d)} (x,y) = \sum_{i=0}^{M-1} \sum_{j=0}^{N-1} F_{i,j} x^iy^j.$$ Let $R(S) = \sum_{i=0}^{M-1} \sum_{j=0}^{N-1} F_{i,j}^2$.
Obviously $R(S) \leq 2^k$, but the following examples show that it can be much smaller. \\
If $S = \{(1,0),(0,1)\}$,then $R(S) = 4$, \\ if $S = \{(1,0), (0,1), (1,1), (1, -1)\}$, then $R(S) = 8$,\\ if $S=\{(1,0), (0,1), (1,1), (1,-1), (2,1), (2,-1), (1,2), (1,-2)\}$, then \\$R(S) = 24$.

\begin{theorem} 
\label{Z3} 
Let $h: A \to  \mathbb{R}$ have integer line sums with respect to $A$. Then there exists a function $f: A \to \mathbb{Z}$ with the same line sums with respect to $A$ satisfying
$$
|\vec{f}-\vec{h}| \leq \frac{1}{2} \sqrt{R(S)(m-M)(n-N)}.
$$
\end{theorem}

In the proof we shall use the following lemma.

\begin{lemma}
\label{lempar}
Let be given a $d$-dimensional parallelepiped $\mathcal{P}$ whose edges have lengths $l_1, l_2, \dots, l_d$. Then the distance from any point of the parallelepiped to the nearest vertex is at most 
$\sqrt{l_1^2+l_2^2+ \dots +l_d^2}/2$.

\end{lemma}

\noindent Note that the bound is the best possible and that it is attained by the centre point of a hyperblock with edge lengths $l_1, l_2, \dots, l_d$. 

\noindent Note further that each face $\mathcal{Q}$ of $\mathcal{P}$ is generated by $d-1$ out of the $d$ edge directions, and that all the vertices of $\mathcal{Q}$ are vertices of $\mathcal{P}$.

\begin{proof}
By induction on $d$. 
%We indicate the directions corresponding to the edges of lengths $l_1, l_2, \dots, l_d$ by $1, 2, \dots, d$, respectively.
%Without loss of generality we may assume that $l_1 \geq l_2 \dots \geq l_d$.
%For $d=1$ the statement is obvious.
%Let $d > 1.$ 
Let $P$ be a point of the parallelepiped $\mathcal{P}$. 
Let $Q$ be the point on the boundary of $\mathcal{P}$ which is nearest to $P$. Let $i$ be the direction which does not occur among the generators of a face $\mathcal{Q}$ to which $Q$ belongs.
Since on the line through $P$ into the direction $i$ the point $P$ is in between two boundary points of $\mathcal{P}$ at distance $l_i$, the distance between $P$ and $Q$ is at most $l_i/2$.

%Project $P$ orthogonally onto the opposite faces of $\mathcal{P}$ generated by the directions $1, \dots, d-1$.
%Call the projections $Q$ and $R$ with the distance between $P$ and $R$ not smaller than the distance between $P$ and $Q$. 
%Since the distance between both faces is at most $l_d$, the distance between $P$ and $Q$ is at most $l_d/2$.
%Let $\mathcal{Q}$ be the face of $\mathcal{P}$ containing $Q$.

%Suppose first that $Q$ is a point on $\mathcal{Q}$. Then, by the induction hypothesis, the distance from $Q$ to the nearest vertex of $\mathcal{Q}$ is at most  $\sqrt{l_1^2 + \dots + l_{d-1}^2}/2$.
%Hence, by Pythagoras' theorem, the distance from $P$ to the nearest vertex of $\mathcal{P}$ is at most $\sqrt{l_1^2 + \dots + l_{d}^2}/2$.

%Suppose now that $Q$ is not a point of $\mathcal{Q}$. Then the segment $PQ$ contains a point $Q_1$ of $\mathcal{P}$.
%Hence the distance between $P$ and $Q_1$ is less than $l_d/2$.
%Let $Q_2$ be the point of $\mathcal{P}$ nearest to $P$. Then the distance from $P$ to $Q_2$ is less than $l_d/2$.
%Further $Q_2$ is on a face $\bf{Q^*}$ of $\mathcal{P}$ generated by the $d$ edge directions except for direction $i$, say.
We claim that $Q$ is the orthogonal projection of $P$ onto the \\$(r-1)$-dimensional hyperplane containing $\mathcal{Q}$.
Let $Q'$ be the orthogonal projection of $P$ onto $\mathcal{Q}$. If $Q \not= Q'$, then the distance between $P$ and $Q'$ is smaller than the distance between $P$ and $Q$.
%If $Q'$ is not equal to $Q_2$ then either $Q_3$ is on $\mathcal{P}$ and is nearer to $P$ than $Q_2$, which is impossible, or
As $Q'$ is not in $\mathcal{P}$, there is a point on the boundary of $\mathcal{P}$ which is in between $P$ and $Q'$. This is impossible.
Thus $Q$ is the orthogonal projection of $P$ onto the face $\mathcal{Q}$ of $\mathcal{P}$ and the distance between $P$ and $Q$ is at most $l_i/2$.
By the induction hypothesis applied to $\mathcal{Q}$, the distance from $Q$ to the nearest vertex of $\mathcal{Q}$ is at most  $\sqrt{l_1^2 + \dots + l_{i-1}^2 + l_{i+1}^2 + \dots + l_d^2}/2$.
Hence, by Pythagoras' theorem, the distance from $P$ to the nearest vertex of $\mathcal{Q}$ is at most $\sqrt{l_1^2 + \dots + l_{d}^2}/2$. 
Since every vertex of $\mathcal{Q}$ is a vertex of $\mathcal{P}$, the lemma follows.
\end{proof}

\begin{proof}[Proof of Theorem \ref{Z3}] 
By Lemma 2.3 there exists a function $g:A \to \mathbb{Z}$ satisfying the same line sums as $h$ with respect to $S$.
It follows from Lemma 2.1 that the 
$(m-M)(n-N)$ vectors in the $mn$-dimensional linear space corresponding to the $m$ by $n$ blocks $g_{(u,v,S)} $  generate a lattice $L$ in the $(m-M)(n-N)$-dimensional subspace of vectors which correspond to the $m$ by $n$ blocks of integers with all the line sums equal to 0. 
%It follows from Lemma 2.2 that there exists an integer solution $g$ with the required line sums.
Thus the set of integer solutions is given by $g + L$. 
This grid generates parallelepipeds with edges of length $\sqrt{R(S)}$ and integer solution vectors as vertices which cover the entire real solution space.
Therefore the real solution vector $\vec{h}$ is in one of these parallelepipeds. 
It follows from Lemma \ref{lempar} that the distance from $\vec{h}$ to its nearest lattice point is at most $\frac12 \sqrt{R(S)(m-M)(n-N)}$. 
Since the nearest lattice point corresponds to an integer solution $\vec{f}$, this proves the theorem.
\end{proof}

If there exists a solution $f$ with all entries in $\{0,1\}$, then this solution is the shortest among all integer solutions, since 
$$ D = \sum_{i=1}^m \sum_{j=1}^n f(i,j) \leq \sum_{i=1}^m \sum_{j=1}^{n} f(i,j)^2 = | \vec{f}|^2$$
and equality holds if and only if $f: A \to \{0,1\}$. If there does not exist such a solution, one may ask for the shortest integer solution.
The following upper bound for the Euclidean length of the shortest integer solution 
is an improvement of Theorem 2 of \cite{ht1} where the corresponding upper bound is $(2^{k-1}+1) \sqrt{mn}$.

\begin{cor} 
\label{Z3min}
Let $f:A \to \mathbb{R}$ and a set of directions $S$ be given such that the line sums of $f$ with respect to $S$ are integers. Let $\vec{f_0}$ be the projection vector with respect to $S$ and the line sums.
Let $g: A \to  \mathbb{Z}$ be the shortest integer solution in the Euclidean sense having the same row and column sums as $g$ has. Then 
$$
|\vec{g}-\vec{f_0}| \leq \sqrt{R(S)(m-M)(n-N)}.
$$
\end{cor}

\begin{proof}
Apply Theorem \ref{Z3} with $\vec{h}=\vec{f_0}$. This yields a solution \\
$g: A \to  \mathbb{Z}$ with $|\vec{g}-\vec{f_0}| \leq \sqrt{R(S)(m-M)(n-N)}$. 
Since the shortest integer solution is the integer solution which is nearest to $\vec{f_0}$, the inequality holds for the shortest solution too.
\end{proof}

Lower bounds for $|\vec{g}-\vec{f_0}| $ may be obtained by adapting the methods from Section 4 and applying the inequality of Cauchy-Schwarz.

\section{The torus case}
Let $n$ be a positive integer and
$$
A:=\{(i,j)\in{\mathbb Z}^2\ :\ 0\leq i < n,\ 0\leq j < n\}.
$$
In what follows, by $(u,v)\in{\mathbb Z}^2$ we shall mean the unique point $(i,j)$ of $A$ with $i\equiv u\pmod{n}$ and $j\equiv v\pmod{n}$. (That is, we identify $A$ with the torus ${\mathbb Z}^2\pmod{n}$.)

A direction $(a,b)\in{\mathbb Z}^2$ is admissible for $A$ if
$0\leq a < n$, $0\leq b < n$ and $\gcd(a,b)=1$.
Let $(a,b)$ be an admissible direction for $A$. Then the torus lines into the direction $(a,b)$ in $A$ are defined by
$$
l_{(a,b),(i,j)}=\{(x,y)\in A\ :\ bx-ay\equiv bi-aj\pmod{n}\}\ \ \ ((i,j)\in A).
$$

The next statement shows that torus lines can be obtained from ordinary lines of ${\mathbb Z}^2$.

\begin{lemma}\label{lem-1} Let $(x,y)\in A$, and $l_{(a,b),(i,j)}$ be a torus line in $A$. Then $(x,y)\in 
l_{(a,b),(i,j)}$ if and only if there exist integers $i_0,j_0$ such that $x\equiv i_0\pmod{n}$, $y\equiv j_0\pmod{n}$ and $bi_0-aj_0=bi-aj$.
\end{lemma}

\begin{proof} Since $x\equiv i_0\pmod{n}$, $y\equiv j_0\pmod{n}$ and $bi_0-aj_0=bi-aj$ imply $bx-ay\equiv bi-aj\pmod{n}$, the 'if' part is obvious.

To prove the 'only if' part, let $(x,y)\in l_{(a,b),(i,j)}$, i.e. $bx-ay\equiv t\pmod{n}$ with $t:=bi-aj$. Then there exists an integer $c$ such that
\begin{equation}
\label{lineq}
bx-ay=cn+t.
\end{equation}
Since $\gcd(a,b)=1$, there are integers $x_0$ and $y_0$ such that $c=ay_0-bx_0$. Substituting this into \eqref{lineq}, we obtain $b(x+nx_0)-a(y+ny_0)=t$. Writing $i_0=x+nx_0$ and $j_0=y+ny_0$, the statement follows.
\end{proof}

\begin{lemma}\label{lem0} Let $l_{(a,b),(i,j)}$ be a torus line in $A$. Then we have
$$
l_{(a,b),(i,j)}=\{(i,j)+t(a,b)\ :\ t=0,1,\dots,n-1\}.
$$
\end{lemma}

\begin{proof} Obviously, $\{(i,j)+t(a,b)\ :\ t=0,1,\dots,n-1\}\subseteq l_{(a,b),(i,j)}$. Let now $(x,y)\in A$ be any point of $l_{(a,b),(i,j)}$. Then by Lemma \ref{lem-1} there are integers $i_0,j_0$ with $i_0\equiv x\pmod{n}$ and $j_0\equiv y\pmod n$ such that we have
\begin{equation}\label{eqeq}
bi_0-aj_0=bi-aj.
\end{equation}
By $\gcd(a,b)=1$ this implies $a\mid i_0-i$ and $b\mid j_0-j$, whence $i_0=i+t_1a$ and $j_0=j+t_2b$ for some integers $t_1,t_2$. If $ab=0$ then we clearly may take $t_1=t_2$. Otherwise, from \eqref{eqeq} we get $bt_1a-at_2b=0$ implying $t_1=t_2$. Thus in any case we have
$$
(x,y)\equiv (i_0,j_0)\equiv (i,j)+t(a,b)\pmod{n}
$$
for some $t\in\{0,1,\dots,n-1\}$. Hence the statement follows.
\end{proof}

\begin{lemma}\label{lem1} Every torus line in $A$ goes through precisely $n$ points of $A$. Further, two distinct torus lines in $A$ into the same direction have no points in common.
\end{lemma}

\begin{proof} Let $l_{(a,b),(i,j)}$ be a torus line in $A$.
Noting that for $0\leq t_1<t_2<n$, the points $(i,j)+t_1(a,b)$ and $(i,j)+t_2(a,b)$ are obviously distinct, the first part of the statement immediately follows from Lemma \ref{lem0}.

To prove the second statement, assume that $(i_1,j_1)+t_1(a,b)=(i_2,j_2)+t_2(a,b)$ with some $(i_1,j_1),(i_2,j_2)\in A$ and $0\leq t_1\leq t_2<n$. Then we have $(i_1,j_1)=(i_2,j_2)+(t_2-t_1)(a,b)$, whence obviously,
$l_{(a,b),(i_1,j_1)}\subseteq l_{(a,b),(i_2,j_2)}$, which by symmetry implies the second statement.
\end{proof}

In view of Lemma \ref{lem1}, for any admissible direction $(a,b)$ for $A$ there exist $n$ points $P_1(a,b),\dots,P_n(a,b)$ in $A$ such that the torus lines with direction $(a,b)$ through these points are parallel, and moreover, each point of $A$ belongs precisely to one of these lines.

If $g:\ A\to {\mathbb R}$ is a function, then the torus line sums of $g$ along an admissible direction $(a,b)$ for $A$ are defined by
$$
L_g((a,b),t):=\sum\limits_{(i,j)\in l_{(a,b),P_t(a,b)}} g(i,j)\ \ \ (t=1,\dots,n).
$$
Let $S$ be an arbitrary set of admissible directions for $A$ and $g$ be any fixed function mapping $A$ into ${\mathbb R}$. Then the functions $f:\ A\to {\mathbb R}$ with
$$
L_f((a,b),t)=L_g((a,b),t)\ \ \ ((a,b)\in S,\ t=1,\dots,n)
$$
clearly can be identified with a linear manifold in ${\mathbb R}^{n^2}$. Write $H_g$ for the set of such functions $f$. Further, write
$$
D_g:=\sum\limits_{t=1}^n L_g((a,b),t)
$$
for the total sum of $g$ (sum of the line sums in a direction $(a,b)$). Observe that $D_g$ is independent of the choice of $(a,b)$.

\begin{lemma}\label{lem2} Let $(i,j),(u,v)\in A$ and $(a,b),(c,d)$ be admissible directions for $A$ with $\gcd(ad-bc,n)=1$. Then the torus lines $l_{(a,b),(i,j)}$ and $l_{(c,d),(u,v)}$ have precisely one point in common.
\end{lemma}

\begin{proof} Clearly, the common points of the torus lines $l_{(a,b),(i,j)}$ and $l_{(c,d),(u,v)}$ belong to the solutions $(t,r)$ of the system of congruences
$$
\begin{cases}
ta-rc\equiv u-i\pmod{n}\\
tb-rd\equiv v-j \pmod{n}
\end{cases}
.
$$
Since by $\gcd(ad-bc,n)=1$ the above system has precisely one solution $(t,r)$ modulo $n$, the statement follows.
\end{proof}

Two admissible directions $(a,b)$ and $(c,d)$ for $A$ are called independent, if $\gcd(ad-bc,n)=1$.

\begin{theorem} Let $S=\{(a_1,b_1),\dots,(a_k,b_k)\}$ be a set of pairwise independent admissible directions for $A$, and let $g:\ A\to {\mathbb R}$. Then the shortest element $f_0\in H_g$ (i.e. the one with smallest Euclidean norm) is given by
$$
f_0(i,j)=\sum\limits_{(a,b)\in S}{\frac{L_g((a,b),t(i,j))}{n}}-{\frac{(k-1)T_g}{n^2}} \ \ \ ((i,j)\in A)
$$
where $t(i,j)$ is that index from $\{1,\dots,n\}$ for which $(i,j)$ belongs to the line $l_{(a,b),P_{t(i,j)}}$.
\end{theorem}

\begin{proof} To prove the statement, we need to check two properties: that $f_0$ is a solution, and that $f_0$ is orthogonal to $H_0$, i.e. to the subspace $H_g$ with $g$ being identically zero.

We start with the first property. Take any torus line
$$
l_{(a,b),P_t(a,b)}=P_t(a,b)+r(a,b)\ \ \ (t\in\{1,\dots,n\},\ r=0,1,\dots,n-1).
$$
Further, let $(c,d)\in S$ with $(c,d)\neq (a,b)$. Then by Lemma \ref{lem2}, each torus line $P_s(c,d)+r(c,d)$ has precisely one point in common with $l_{(a,b),P_t(a,b)}$. 
By the definition of $f_0$ this yields, for $(a,b) \in S$,
$$
\sum\limits_{r=1}^n f_0(P_t(a,b)+r(a,b))=
L_g((a,b),t)+{\frac{(k-1)T_g}{n}}-{\frac{(k-1)T_g}{n}},
$$
which shows that $f_0\in H_g$.

To prove the second property, let $h\in H_0$ be arbitrary. Then we clearly have
$$
\sum\limits_{r=0}^{n-1} h(P_t(a,b)+r(a,b))=0\ \ \ ((a,b)\in S,\ t=1,\dots,n).
$$
Thus
$$
\sum\limits_{(i,j)\in A} f_0(i,j)h(i,j)=
\sum\limits_{(i,j)\in A} \left(\sum\limits_{(a,b)\in S}{\frac{L_g((a,b),t(i,j))}{n}}-{\frac{(k-1)T_g}{n^2}}\right) h(i,j)=
$$
%$$
%=\sum\limits_{(a,b)\in S} \sum\limits_{s=1}^n \sum\limits_{r=0}^{n-1}
%\left({\frac{L_g((a,b),s)}{n}}-{\frac{(k-1)T_g}{n^2}}\right) h(P_s(a,b)+r(a,b))=
%$$
$$
=\sum\limits_{(a,b)\in S} \sum\limits_{s=1}^n \left({\frac{L_g((a,b),s)}{n}}-{\frac{(k-1)T_g}{n^2}}\right) \sum\limits_{r=0}^{n-1}h(P_s(a,b)+r(a,b))=0
$$
and the theorem follows.
\end{proof}

\section{The continuous version in the simple case}
%Finally we prove an analogue of Theorem \ref{thm2}.
Let $m$ and $n$ be positive real numbers and let $A$ be a
Lebesgue-measurable subset of $T=[0,m]\times [0,n]$. Write
$f_A(x,y)$ for the characteristic function of $A$ inside $T$.

For any bounded Lebesgue-measurable function $f:\ T\to {\mathbb R}$, let $c_f(x):\ [0,m]\to {\mathbb R}$ and $r_f(y):\ [0,n]\to {\mathbb R}$ be the column integrals and row integrals of
$f(x,y)$, that is
$$
c_f(x)=\int\limits_0^n f(x,y) dy\ \ \ (x\in [0,m])
$$
and
$$
r_f(y)=\int\limits_0^m f(x,y) dx\ \ \ (y\in [0,n])),
$$
respectively, where integration is always meant in the Lebesgue sense. Since $f(x,y)$ is bounded, by the theorem of Fubini we know that these functions exist. Note that the same is true for $c_{f_A}(x)$ and $r_{f_A}(y)$.

Let $L$ denote the set of bounded Lebesgue-integrable functions $T\to {\mathbb R}$ having column integrals $c_{f_A}(x)$ $(x\in [0,m])$ and row integrals $r_{f_A}(y)$ $(y\in [0,n])$. Further, write $H$ for the set of bounded Lebesgue-integrable functions $T\to {\mathbb R}$ having vanishing row integrals and column integrals. Observe that $H$ is a closed linear subspace of the linear space ${\mathcal L}$ of bounded integrable functions $T\to {\mathbb R}$. Further, for any $g_1,g_2\in L$ we obviously have $g_1-g_2\in H$. In other words, $L=g+H$ with any $g\in L$.

We recall that the inner product in ${\mathcal L}$ is given by
$$
\langle f(x,y),g(x,y)\rangle =\underset{T}{\int\int} f(x,y)g(x,y) dxdy
$$
for $f,g\in {\mathcal L}$. The following theorem describes the shortest element in $L$, with respect the usual norm
$$
||f(x,y)||=\sqrt{\langle f,f\rangle}= \left(\underset{T}{\int\int} f^2(x,y) dxdy\right)^{1/2}
$$
for $f\in {\mathcal L}$.

\begin{theorem}
\label{thm6} The shortest element in $L$ exists, and is given by
$$
f_0(x,y)={\frac{c_{f_A}(x)}{n}}+{\frac{r_{f_A}(y)}{m}}
-{\frac{\lambda(A)}{mn}}\ \ \ \ ((x,y)\in T),
$$
where $\lambda(A)$ is the Lebesgue-measure of $A$.
\end{theorem}

\begin{proof}
Since $H$ is a closed linear subspace of ${\mathcal L}$ and $L$ is just a shift of $H$, $L$ has a shortest element $f_0(x,y)$ indeed. This $f_0(x,y)$ is uniquely determined by the following two properties:
\begin{itemize}
\item $f_0(x,y)$ has column integrals $c_{f_A}(x)$ $(x\in [0,m])$ and row integrals $r_{f_A}(y)$ $(y\in [0,n])$,
\item $f_0(x,y)$ is orthogonal to $H$, i.e. $\langle f_0(x,y),h(x,y)\rangle=0$ for every $h(x,y)\in H$.
\end{itemize} 

We prove that the choice for $f_0(x,y)$ in the statement meets these requirements. To prove the first property, observe that
$$
\int\limits_0^n r_{f_A}(y) dy=\underset{T}{\int\int} f_A(x,y) dxdy=\lambda(A).
$$
Thus for any $x\in [0,m]$ we have
$$
\int\limits_0^n \left({\frac{c_{f_A}(x)}{n}}+{\frac{r_{f_A}(y)}{m}}-
{\frac{\lambda(A)}{mn}}\right) dy=
c_{f_A}(x)+{\frac{\lambda(A)}{m}}-n{\frac{\lambda(A)}{mn}}=
c_{f_A}(x).
$$
Similarly, for any $y\in [0,n]$
$$
\int\limits_0^m \left({\frac{c_{f_A}(x)}{n}}+{\frac{r_{f_A}(y)}{m}}-
{\frac{\lambda(A)}{mn}}\right)dx=
{\frac{\lambda(A)}{n}}+r_{f_A}(y)-m{\frac{\lambda(A)}{mn}}=
r_{f_A}(y),
$$
which proves the first property.

In order to check the second property, take an arbitrary $h(x,y)\in H$. For any $x\in [0,m]$ and $y\in [0,n]$
we have
$$
c_h(x)=\int\limits_0^n h(x,y) dy=0 ~~ {\rm and} ~~
r_h(y)=\int\limits_0^m h(x,y) dx=0.
$$
Then
$$
\left\langle h(x,y), {\frac{c_{f_A}(x)}{n}}+{\frac{r_{f_A}(y)}{m}}-
{\frac{\lambda(A)}{mn}}\right\rangle=
$$
$$
=\int\limits_0^n\int\limits_0^m h(x,y)\left({\frac{c_{f_A}(x)}{n}}+{\frac{r_{f_A}(y)}{m}}-
{\frac{\lambda(A)}{mn}}\right)dxdy=
$$
$$
\int\limits_0^m{\frac{c_{f_A}(x)}{n}}
\left(\int\limits_0^n h(x,y) dy\right)dx+
%$$
%$$
\int\limits_0^n
\left({\frac{r_{f_A}(y)}{m}}-{\frac{\lambda(A)}{mn}}\right)
\left(\int\limits_0^m h(x,y)dx\right)dy
$$
%$$
%=\int\limits_0^m{\frac{c_{f_A}(x)}{n}}c_h(x)dx
%+\int\limits_0^n \left({\frac{r_{f_A}(y)}{m}}-{\frac{\lambda(A)}{mn}}\right)
%r_h(y)dy=0.
%$$
and the inner integrals are 0. This proves the second property, and the theorem follows.
\end{proof}

\end{document}